\newcommand{\rr}{\mathbb{R}}
\newcommand{\p}{\partial}
\newcommand{\tn}{|\!|\!|}

\documentclass[12pt]{amsart}
\usepackage{amsmath}
\usepackage{amsthm}
\usepackage{amssymb}
\usepackage{eucal}
\usepackage[left=1.1in,top=2in,right=1.1in, bottom=1.5in]{geometry}

\begin{document}

\newtheorem{theorem}{Theorem}[section]
\newtheorem{lemma}[theorem]{Lemma}
\newtheorem{proposition}[theorem]{Proposition}
\newtheorem*{lem1}{Lemma}

\numberwithin{equation}{section}

\title[Schr\"odinger Equations]
{On the decay properties of solutions to a class of Schr\"odinger
equations}

\author{L. Dawson}
\address{Department of Mathematics,
University of California
Santa Barbara, CA 93106
}
\email{ldawson@math.ucsb.edu}

\author{H. McGahagan}
\address{Department of Mathematics,
University of California
Santa Barbara, CA 93106
}
\email{helena@math.ucsb.edu}

\author{G. Ponce}
\address{Department of Mathematics,
University of California
Santa Barbara, CA 93106
}
\email{ponce@math.ucsb.edu}

\thanks{The first and third authors were supported by 
NSF grants, and the second author was supported by an NSF
postdoctoral fellowship.}

\subjclass[2000]{Primary 35J10; Secondary 35B65}

\date{}

\begin{abstract}
We construct a local in time, exponentially decaying solution of 
the one-dimensional 
variable coefficient Schr\"odinger equation by solving a 
nonstandard boundary value problem.  A main ingredient in the 
proof is a new commutator estimate involving the projections
$P_\pm$ onto the positive and negative frequencies. 
\end{abstract}

\maketitle

\section{Introduction}\label{S: intro}

In \cite{Ka}, T. Kato showed that the semigroup
$\{e^{-t\partial_x^3}\,:\,t\geq 0\}$ in the space $L^2(e^{2 \beta
x}dx)$ with $\beta
>0$ is formally equivalent to the semigroup
$e^{-t(\partial_x-\beta)^3}$ in $L^2(\mathbb{R})$.  Among the
immediate consequences of this result is that if $u\in C([0,T] :
H^1(\mathbb{R}))$ is a strong solution of the $k$-generalized
Korteweg de Vries (KdV) equation,
\begin{equation}\label{E: kdv}
\partial_tu+\partial_x^3u+u^k\,\partial_xu=0, \quad k=1,2,\ldots
\end{equation}
with data $u_0\in L^2(e^{2 \beta x}dx)$, then $u\in C([0,T] :
L^2(e^{2\beta x}dx)) \cap C^{\infty}(\mathbb{R}\times (0,T])$. In other words,
the solution $u=u(x,t)$ satisfies the persistence property $
e^{\beta x}u \in C([0, T] : L^2(\mathbb{R}))$ and a ``parabolic"
regularization, $u\in C^{\infty}(\mathbb{R}\times(0, T]).$

Since results for solutions of  the $k$-generalized KdV equation
and  Schr\"{o}dinger equations of the type
\begin{equation}
\label{E:1} (a)\;\;\partial_tu -i\Delta u = f(|u|)u,\quad \quad
(b)\;\;\partial_tu -i(\Delta u+W(x,t)u)=F(x,t),
\end{equation}
run parallel -- for instance, solutions of both satisfy 
Strichartz estimates, local smoothing
effects of the Kato type, and persistence properties in
$H^s(\mathbb{R})$, the weighted spaces $H^s(\mathbb{R})\cap
L^2(|x|^k)$, and the Schwartz space -- one may ask what the
equivalent result to that described above for the KdV equation is in
the case of Schr\"{o}dinger equations. One first notices that even
for the free Schr\"odinger group $\{ e^{it\Delta}\,:\,t\in
\mathbb{R}\}$, both of the above properties fail: assuming we
are in $\mathbb{R}^1$ $(\Delta = \partial_x^2)$ for simplicity,
we can construct initial data
$u_0\in L^2(\mathbb{R})\cap L^2(e^{2\beta x}dx)$ such
that $e^{it\p_x^2}u_0\notin L^2(e^{2\beta x}dx)\cup
C^{\infty}(\mathbb{R})$ for any $t>0$.

Roughly, the  difficulty lies in the fact that if
$u(x,t)=e^{it\p_x^2}u_0(x)$, then $v(x,t):=e^{\beta x}
 u(x,t)$ formally solves the equation
 \begin{equation}
 \label{111}
 \p_t v - i(\p_x-\beta)^2 v=\p_t v-i\p_x^2 v + 2 i \beta \p_x v -i\beta^2 v=0,
 \end{equation}
 whose associated initial value problem (IVP) is ill-posed in
$L^2(\mathbb{R})$. However, the operator $2i\beta
\partial_x$, whose symbol is  $-2\beta \xi$, introduces a
parabolic structure in the negative frequency for positive time
and in the positive frequency for negative time. Thus, to find
$L^2$-solutions of equation \eqref{111} in the time interval
$[0,T]$, one needs to consider a ``boundary value problem" 
for \eqref{111} where
 \begin{equation}
 \label{data}
 \begin{aligned}
  &v_{-}(x,0)=P_{-}v(x,0) := 
(\chi_{(-\infty,0)}(\xi)\hat v(\xi,0))^\lor (x),\\
&v_{+}(x,T)=P_{+}v(x,T) := 
(\chi_{(0,\infty)}(\xi)\hat v(\xi,T))^\lor (x)
\end{aligned}
\end{equation}
are prescribed. In this case, 
one finds the solution
\begin{equation}
 \label{112}
 v(x,t)= e^{ t(i\p_x^2-2\beta D_x +i\beta^2)}v_{-}(x,0) +
  e^{ -(T-t)(i\p_x^2+2\beta D_x +i\beta^2)}v_{+}(x,T),
 \end{equation}
with $D_x h(x):=(-\p_x^2)^{1/2}h(x)=(c|\xi| \hat
h(\xi))^{\lor}(x)$.  Then,
\begin{equation}
 \label{freeest}
 \sup_{[0,T]}\|v(t)\|_2\leq c(\|v_{+}(x,T)\|_2+\|v_{-}(x,0)\|_2),
 \end{equation}
 $c$ independent of $\beta>0$ and $T$,
and $v\in C^{\infty}(\mathbb R\times (0,T)).$
We observe 
that
in formula \eqref{112}, the positive and negative
frequencies do not interact and, also, 
that  $u(x,t):=e^{-\beta x} v(x,t)$ is not necessarily an
$L^2$-solution of the free Schr\"odinger equation.

The following estimate established in \cite{KePoVe2} of the type
described in \eqref{freeest} 
for a linear Schr\"odinger equation with
lower order variable coefficients  \eqref{E:1}\rm{ (b)} was a key
step in the proof of the unique continuation results obtained in
\cite{KePoVe2} and \cite{EsKePoVe}.

\begin{lem1}\cite{KePoVe2} There exists $\epsilon>0$ such that if $\, W:\Bbb R^n\times [0,T]\to\Bbb C$ satisfies 
$
\|W\|_{L^1_tL^{\infty}_x}\leq \epsilon
$
and $u\in C([0,T]:L^2_x(\Bbb R^n))$ is a strong solution of the equation \eqref{E:1}{\rm{ (b)}}
with
\begin{equation}
\label{hyp} u_0=u(\cdot,0),\,\,u_T\equiv u(\cdot,T)\in
L^2(e^{2\beta x_1}dx),\quad F\in L^1([0,T]:L^2_x(e^{2\beta
x_1}dx))
\end{equation}
for some $\beta\in\Bbb R$, then there exists $c$ independent of 
$\beta$ such that
$$
\sup_{0\leq t\leq T}\| e^{\beta x_1} u(\cdot,t)\|_2 \leq c
\Big(\|e^{\beta x_1} u_0\|_2 + \|e^{\beta x_1} u_T\|_2 +\int_0^T
\|e^{\beta x_1} F(\cdot, t)\|_2 dt\Big).
$$
\end{lem1}

Notice that in the above result, one assumes the existence of a reference
solution $u(x,t)$ of equation \eqref{E:1} (b) and shows that under
hypothesis \eqref{hyp}, exponential decay in the
time interval $[0,T]$ is preserved.

The $L^2$-well-posedness of the IVP associated to the equation
\begin{equation}
\label{123}
\partial_t w = i\Delta w + b(x)\cdot \nabla_x w + f(x,t),
\end{equation}
has been extensively studied. 
In particular, S. Mizohata \cite{Mi} gives the following 
necessary condition for the IVP associated to \eqref{123} 
to be well-posed in $L^2(\rr^n)$:
\begin{equation}
\label{miz}
\sup_{x\in\mathbb R^n, \omega \in \mathbb S^{n-1}, R>0}\;
|\;\text{Im} \;\int_0^{R}
b(x+r\omega)\cdot\omega dr| < \infty.
\end{equation}
The gain of regularity of solutions to the variable coefficient 
Schr\"odinger equation
$$
\partial_tu -i\p_{x_j}(a_{jk}(x)\p_{x_k}u)+W(x)u= 0
$$
as a consequence of its dispersive character and the decay
assumptions on the data has 
also 
been studied in several works; see
\cite{CrKaSt}, \cite{DS}, and references therein.

In this note, we shall combine the above ideas with some new commutator
estimates to construct an exponentially decaying solution to the
one-dimensional variable coefficient Schr\"odinger equation
\begin{equation}
\label{E2}
\partial_tu= i(\partial_x(a(x,t)\partial_xu)+W(x,t)u).
\end{equation}
More precisely, we are interested in a solution $u \in C([0,T] :
L^2(\mathbb{R}) \cap L^2(e^{2 \beta x}dx))$.

To ensure that we construct $u \in L^2(\rr)$, we will need to refer
to the following function $\varphi_\beta(x)$:
for $\beta>0$ we denote by $ \varphi(x)=\varphi_{\beta}(x)$ a
$C^4(\mathbb R)$ function such that $\varphi(x)=1$
if $x\leq 0$, $\varphi(x)=e^{\beta x}$ if $x\geq10\beta$, and 
$\varphi(x)$ is strictly increasing on $(0,10 \beta)$. 

\begin{theorem}\label{T: T1}
Let  $a:\mathbb R\times \mathbb R^+\to\mathbb R$ be such that
\begin{equation}
\label{hyp2}
\begin{aligned}
& a\in C^2(\mathbb R\times\mathbb R^+)\cap L^1_t(\mathbb R^+: L^{\infty}_x(\mathbb R)),\;\;\langle x\rangle \p_x^j a\in L^1_t(\mathbb R^+: L^{\infty}_x(\mathbb R)),\;\;j=1, 2,\\
&a(x,t) \geq \lambda \geq 0,\;\; \forall\;(x,t)\in \mathbb
{R}\times \mathbb R^+.
\end{aligned}
\end{equation}
Let $W:\mathbb R\times \mathbb R^+\to\mathbb C$
be such that
\begin{equation}
\label{hyp4}
W\in L^1_t(\mathbb R^+: L^{\infty}_x(\mathbb R)).
\end{equation}
  Then given $(f, g) \in P_{-}L^2(\mathbb{R})\times P_{+}L^2(\mathbb{R})$, 
there exists  $T = T(\beta;\tn a\tn_1;
\| W \|_{L^1_t L^\infty_x}) > 0$
such that \eqref{E2} has 
a unique solution  
$u\in C([0,T]:L^2(\mathbb {R}))$ with 
$e^{\beta x} u\in C([0,T] : L^2(\mathbb {R}))$and
with 
$P_{-}(\varphi(x)  u(x,0))=f(x)$ and  
$P_{+}(\varphi(x) u(x,T))=g(x)$. 

If in addition $a, W\in C^{\infty}(\mathbb
R\times\mathbb R^+)$; $\lambda > 0$ with 
\[
\beta \lambda \geq c(
\|\langle x\rangle \p_x a\|_{L^\infty(\rr \times \rr^+)} + 
\|\langle x\rangle \p^2_x a\|_{L^\infty(\rr \times \rr^+)});
\] 
and $\p_t^k\p_x^j a,  \,\p_t^k\p_x^j W\in
L^{\infty}(\mathbb R\times\mathbb R^+) $ for any $k, j\in\mathbb
Z^+$, then $\; u\in C^{\infty}(\mathbb R\times (0,T))$.
\end{theorem}

We use the notation $\langle x \rangle :=
(1+|x|^2)^{1/2}$.  Also, $\tn a \tn _1$ denotes the sum of the
$L_t^1L_x^\infty$-norms of the expressions involving the function $a$
described in \eqref{hyp2}:
\begin{equation*}
\tn a \tn _1 := \|a\|_{L^1_tL_x^{\infty}}+ \sum_{j=1}^{2} \|
\langle x \rangle \partial_x^ja\|_{L_t^1L^{\infty}_x}.
\end{equation*}

Under the assumptions of Theorem \ref{T: T1}, we do not know if the
dependence on the parameter $\beta$ of the time interval $[0,T]$
can removed as was done in \cite{KePoVe2}.
Also, here we shall restrict ourselves to the one-dimensional case.

To prove Theorem \ref{T: T1}, we consider a system  describing the
time evolution of the projection of the weighted function 
$v\!:=\!\varphi u$ into
the positive and negative frequencies. Since our equation has 
variable coefficients, this becomes a coupled system. 
It will be essential in 
our arguments that the coupled terms are,
roughly speaking, of  ``order zero."  We will show this using
commutator estimates such as the following:
for all $p\in(1,\infty),$ $l,\,m\in\mathbb  Z^+$ 
there exists $c=c(p;l;m)>0$ such that
 \begin{equation}
\label{77}
 \| \p_x^l[P_{+};\,a]\p_x^m f\|_p\leq c \|\p_x^{l+m} a\|_{\infty} \|f\|_p.
\end{equation}
Clearly, the inequality \eqref{77} holds with $P_{-}$ or $H$, the
Hilbert transform, in place of $P_{+}$.  In the
case $l+m=1$,
\eqref{77} is Calder\'on's first commutator estimate
\cite{Ca}. A related version of estimate \eqref{77} was obtained
in \cite{MoRi} for general positive  derivatives, but did not
involve the $L^{\infty}$-norm.

\section{Proof of Theorem 1.1}\label{S: calc}

Consider the equation
\begin{equation}
\label{2.1}
\partial_tu= i(\partial_x(a(x,t)\partial_xu)+W(x,t)u).
\end{equation}
We wish to construct a solution $u \in L^2((1+e^{2 \beta x})dx)$ 
for a fixed $\beta>0$.
Recall the definition of the function 
$\varphi(x)=\varphi_\beta(x)$, and define
$\phi(x) := \varphi'(x)/\varphi(x)$. Notice that 
$\phi(x) = \beta \chi_{\rr^+}(x)$ except on 
the interval $0 < x < 10 \beta$ and that  
$\| \phi \|_\infty = \beta$. 

Let $v(x,t) := \varphi(x) u(x,t)$. Then,
multiplying \eqref{2.1} by $\varphi(x)$ and using the
fact that $[\varphi; \p_x] = - \phi \varphi$, we have that
\begin{equation}
\label{2.2}
\begin{aligned}
\partial_xv &= i((\p_x-\phi(x))(a(x,t)(\p_x-\phi(x))v)+W(x,t)v)\\
&=i\partial_x(a\partial_xv)-2ia \phi \partial_xv +i(
(\phi^2-\p_x\phi)a-\phi\partial_x a)v+iWv.
\end{aligned}
\end{equation}
We will construct a solution $v \in L^2(\rr)$ 
of \eqref{2.2}.  This suffices since the definition of 
$\varphi$ then guarantees that $u$ defined by
$u(x) = v(x)$ on $x \leq 0$ and $u(x) = \varphi^{-1}(x) v(x)$
on $x > 0$ will be in $L^2((1+e^{2 \beta x})dx)$, and 
$u$ will solve \eqref{2.1}.

Applying the projection operators $P_\pm$  
to equation \eqref{2.2}, 
we obtain
\begin{align*}
\partial_t v_{\pm} 
&=i\partial_x(a\partial_xv_{\pm})-2i\phi a\partial_xv_{\pm}+
P_{\pm}(i((\phi^2-\p_x\phi) a -\phi \partial_x a
))v)+P_{\pm}(iWv)\\
&\quad+i\partial_x( [P_{\pm};a]\partial_xv) -2i[P_{\pm};
a\phi]\partial_xv,
\end{align*}
where $v_\pm := P_\pm v$.  We 
can rewrite this as the following coupled system:
\begin{equation}
\label{system}
\aligned
\partial_t v_{+} &= i\partial_x(a\partial_xv_{+})
- 2 i a\phi\partial_xv_{+} + \Lambda_{+}(v_{+},v_{-})
\\
\partial_t v_{-} &= i\partial_x(a\partial_xv_{-})
- 2 i a\phi\partial_xv_{-} + \Lambda_{-}(v_{+},v_{-}),
\endaligned
\end{equation}
where
$$
\aligned
\Lambda_{\pm}(v_{+},v_{-})
&:=
P_{\pm}(i((\phi^2-\p_x\phi) a - \phi \partial_x  a)(v_+ + v_-))
+P_{\pm}(i W (v_+ + v_-))
\\ &\quad
+i \partial_x([P_{\pm};a]\partial_x (v_+ + v_-) 
- 2 i[P_{\pm}; a \phi]\partial_x (v_+ + v_-)).
\endaligned
$$
Notice that once we construct functions 
$v_{+}$ and $v_{-}$ that solve this system, 
$v = v_{+} + v_{-}$ will be the desired solution of \eqref{2.2}.

Taking the $L^2$ norm 
of $\Lambda_{\pm}$
and applying Lemma \ref{L: L3}, 
it follows that $\Lambda_{\pm}$ can be written
as a sum of linear operators in $(v_{+},v_{-})$ of ``order zero": 
\begin{equation}
\begin{aligned}
\|\Lambda_{\pm}(v_{+}, v_{-})\|_{2}
&\leq 
c \, \Big( \| (\phi^2 - \p_x \phi)a - \phi \p_x a  \|_\infty + 
	\| W \|_\infty 
+ \|\p_x^2 a \|_\infty 	
\\
&\quad\quad\quad+ \| \p_x (a \phi) \|_\infty \Big)\, 
\|v_{+} + v_{-}\|_2
\leq
K(t)\,(\|v_{+}\|_2 + \|v_{-}\|_2),
\end{aligned}
\end{equation}
with
$$
K(t) := c\,\Big(\,\sum_{j=0}^2
\beta^j\,\|\p_x^{2-j}a(t)\|_{\infty}+\|a(t)\|_{\infty}
+\|W(t)\|_{\infty}\Big).
$$

To prove the existence of a solution $(v_+, v_-) \in L^2$ to 
\eqref{system}, we will 
establish {\it{a priori}} estimates and local existence for a 
related uncoupled system, and then find $(v_+, v_-)$
as a limit of these solutions.

First, we fix the time interval on which we will 
solve the equation.  Define
$$
c_{a,\beta}(t) := c \Big(\|a(t)\|_{\infty}+(1+\beta)\|
\langle x \rangle \p_xa(t)\|_{\infty}+\beta\| \langle x \rangle
\p^2_xa(t)\|_{\infty}\Big),
$$
and let
$\,T = T(\beta; \tn a \tn_1; \|W\|_{L^1_t L^\infty_x})>0\,$ 
be such that
\begin{equation}
\label{TT} e^{4 \int_0^T c_{a,\beta}(t)dt}\,\leq 2/3\quad
\text{and}\quad \int_0^T K(t)dt\leq1/8.
\end{equation}
These inequalities must hold for some $T > 0$ 
by hypotheses \eqref{hyp2} and 
\eqref{hyp4}. 
Also, we define the norm
$
\tn v \tn _T := \sup_{[0,T]} \|v_+(t)\|_2 +\sup_{[0,T]} \|v_-(t)\|_2,
$
and letting $\delta := \| v_+(T)\|_2 + \|v_{-}(0)\|_2$, 
we define the space 
\begin{equation}
\label{space}
\mathrm{X}_T := \{v: \mathbb{R}\times [0,T] \to \mathbb{C}: 
\tn v \tn _T \leq 
4 \delta \}.
\end{equation}

Next, 
using standard energy estimates, we obtain {\it{a priori}} 
bounds for the solutions of both of the following 
(uncoupled) equations on $\rr \times [0,T]$:
\begin{eqnarray}
\label{eq11}
\p_t v_{+} 
&=& 
i\partial_x(a \partial_xv_{+})- 2i a \phi \partial_xv_{+}
+ F_{+}(x,t)
\\
\label{eq22}
\p_t v_{-} 
&=&
i\partial_x(a\partial_xv_{-})-2i a \phi \partial_xv_{-} +F_{-}(x,t),
\end{eqnarray}
with functions $F_{\pm} \in L^1_t(\rr^+:L^2_x(\rr))$.
Multiplying \eqref{eq22} by $\overline{v_{-}}$, integrating in
the $x$-variable, and taking the real part, 
we have that
\[
\frac{1}{2} \frac{d\;}{dt} \|v_{-}(t)\|_2^2
= Re \left\{ - 2 i \int a \phi \p_x v_{-} \overline{v_{-}}dx 
+ \int F_{-}(x,t) \overline{v_{-}}dx \right\}.
\]
Using the definition of $D_x^\alpha$ 
and the fact that $\widehat v_-$ 
is supported on $\rr^-$, we compute
$$
-2i  \int a \phi \p_x v_{-} \overline{v_{-}} dx
=-2  \int a \phi | D_x^{1/2} v_{-} |^2 dx 
-2 \int (D_x^{1/2} ([D_x^{1/2};a\phi ] \overline{v_{-}}) v_{-} dx
;
$$ 
therefore, 
\begin{equation}
\label{estimate}
\begin{aligned}
\frac{d\;}{dt}\|v_{-}(t)\|_2^2
+ 4 &\int  a \phi |D_x^{1/2} v_{-}|^2 
\leq
4 c_{a,\beta}(t)\,\|v_{-}(t)\|_2^2 +
2 \|F_{-}(t)\|_2 \|v_{-}(t)\|_2,
\end{aligned}
\end{equation}
where the final inequality follows from
combining the estimate
from Proposition \ref{P: P3} in the appendix and 
the Gagliardo-Nirenberg inequality  to see that
$$
\begin{aligned}
\|D_x^{1/2} [D_x^{1/2};a\phi] v_{-}\|_2
&\leq
c \|J^\delta \p_x (a \phi)\|_{q} \|v_{-}\|_2
\leq 
c \|\p_x (a\phi) \|_{q}^{1-\delta} 
\|J \p_x (a\phi)\|_{q}^\delta \| v_{-}\|_2
\\
&\leq
c \Big( \| \p_x (a \phi) \|_{q} + \| \p_x^2 (a \phi)\|_{q} \Big) 
\| v_{-}(t) \|_2
\leq 
c_{a,\beta}(t)\|v_{-}(t)\|_2,
\end{aligned}
$$
where we take $q < \infty$ and $0 < \delta < 1$ such that
both $\delta > 1/q$ and 
$\delta > 1 - 1/q$, and also $q$ large enough that 
$\| \langle x \rangle \|_q < \infty$.  
Bounding $\frac{d \;}{dt} \| v_-(t) \|_2$ from (2.9), we find 
that 
$\|v_{-}(t) \|_2 \leq 
\Big( \| v_{-}(0)\|_2 + \int_0^T \| F_{-} \|_2\Big) 
e^{2 \int_0^T c_{a, \beta}(\tau) d\tau}$ for all $t \in [0, T]$.
Putting this back into (2.9) in order to bound
$\int_0^T \int a \phi | D_x^{1/2} v_-|^2 \;dxdt$, 
we obtain the estimate
\begin{equation}
\label{estimate1}
\begin{aligned}
\sup_{t \in [0, T]} \|v_-(t)\|_2 + 
2 \Big( & 
\int_0^T \int a(x,t) \phi(x) |D_x^{1/2}v_{-}|^2\; dx dt \Big)^{1/2}
\\
&\leq 
3 \Big(\|v_{-}(0)\|_2+\int_0^T\|F_{-}(t)\|_2
\;dt\Big)\,e^{4 \int_0^Tc_{a,\beta}(\tau)d\tau}.
\end{aligned}
\end{equation}
A similar argument applied to the equation for $v_+$ \eqref{eq11}
shows that
\[
\label{preestimate2}
\frac{d\;}{dt}\|v_{+}(t)\|_2^2-  
4 \int  a\phi |D_x^{1/2} v_{+}|^2\;dx\\
\geq - 4 \,c_{a,\beta}(t) \|v_{+}(t)\|_2^2
-2 \|F_{+}(t)\|_2 \|v_{+}(t)\|_2.
\]
Integrating from $t$ to $T$, we estimate 
$\|v_+(t)\|_2 \leq  
\Big(\| v_+(T) \|_2 + 
\int_0^T \| F_+ \|_2\Big) e^{2 \int_0^T c_{a, \beta}}$, and 
then, it follows that
\begin{equation}
\label{estimate2}
\begin{aligned}
\sup_{t\in[0,T]}\|v_{+}(t)\|_2 +  2 \Big(&\int_0^T\int a(x,t) \phi(x)
|D_x^{1/2}v_{+}|^2\;dx dt\Big)^{1/2} 
\\
&\leq 3
\Big(\|v_{+}(T)\|_2+\int_0^T\|F_{+}(t)\|_2 \;dt\Big)\,e^{4 \int_0^T
c_{a,\beta}(\tau)d\tau}.
\end{aligned}
\end{equation}

To establish the first part of Theorem \ref{T: T1}, the
existence and uniqueness of a solution of \eqref{2.1}, we apply
the contraction principle in the space $X_T$ \eqref{space} with
$(v^{m}_{+},v^{m}_{-})$ for $m \in \mathbb{N}$ the iteratively defined 
solution of the system
  \begin{equation}
  \label{system2}
  \begin{cases}
\begin{aligned}
&\partial_t v_{+}^{m} =i\partial_x(a\partial_xv_{+}^{m})-2i a\phi \partial_xv_{+}^{m} +\Lambda_{+}(v^{m-1}_{+},v^{m-1}_{-}),\\
&\partial_t v_{-}^{m} =i\partial_x(a\partial_xv_{-}^{m})-2i a\phi \partial_xv_{-}^{m} + \Lambda_{-}(v^{m-1}_{+},v^{m-1}_{-}),\\
&v_{+}^{m}(x,T)=g(x),\quad v_{-}^{m}(x,0)=f(x),
\end{aligned}
\end{cases}
\end{equation}
where $v_+^0 = v_-^0 :=0$.
The above equations are of the form 
\eqref{eq11} and \eqref{eq22}, and 
the existence of solutions 
in $C([0, T]: L^2(\rr))$ will be proven below. 
Letting $\|v(t)\|_2 := \| v_+(t)\|_2 + \| v_-(t)\|_2$, 
we have, 
from the 
energy estimates \eqref{estimate1} and \eqref{estimate2},
that
\begin{equation}
\label{energydifference}
  \aligned
  \sup_{t\in[0,T]}\|v^{m+1}(t)\|_2
  &\leq
  3\Big(\delta + 2 \sup_{t\in[0,T]}\|(v^{m}(t)\|_2 
\int_0^T K(t)\; dt\Big)
  \,e^{4 \int_0^Tc_{a,\beta}(\tau)d\tau}\\
  \endaligned
\end{equation}
for $m \in \mathbb{N}$.
From our choice of $T$ in \eqref{TT}, 
$\sup_{t \in [0, T]} \|v^1(t)\|_2 
\leq 3 \delta e^{4 \int_0^T c_{a, \beta}} \leq 2\delta,$ 
and if we assume
$\sup_{t \in [0, T]} \| v^m(t) \|_2 \leq 4\delta$, 
then the energy estimate \eqref{energydifference} 
yields \[
\sup_{t \in [0, T]} \|(v^{m+1}(t)\|_2 
\leq 
3 (\delta+ 2 (4\delta) (1/8)) 2/3 = 4\delta.
\] 

Repeating the derivation of the energy estimates for the 
equations for the differences $v_+^{m+1} - v_+^{m}$ and
$v_-^{m+1} - v_-^m$ and using \eqref{TT} yields the estimate
$$
 \aligned
  \sup_{t\in[0,T]}\|(v^{m+1}-v^m)(t)\|_2
 & \leq 
\frac{1}{2}\,\sup_{t\in[0,T]}\|(v^{m}-v^{m-1})(t)\|_2 .
  \endaligned
  $$
Therefore, by the contraction principle there exists a
unique solution $(v_{+},v_{-})\in C([0,T]:L^2(\mathbb R))$ of the
system \eqref{system} 
(which is realized in $C([0,T]:H^{-2}(\mathbb R))$ 
with data $v_{+}(x,T)=g(x)$ and $v_{-}(x,0)=f(x)$.

To complete the above argument, we shall use the artificial
viscosity method to prove the existence of solutions of
\eqref{eq22} with initial data specified at $0$ (similarly, we can
prove the existence of solutions to \eqref{eq11} 
with data specified at time $T$). 
Thus, we consider the family of equations
\begin{equation}\label{vis}
\partial_tv_{-}^\epsilon 
= -\epsilon \partial_x^4 v_{-}^\epsilon 
+ i \partial_x(a\partial_x v_{-}^\epsilon)  
- 2i a\phi  \p_xv_{-}^\epsilon + F_{-}\\
= -\epsilon \partial_x^4 v_{-}^\epsilon + \Phi (v_{-}^\epsilon), 
\quad t>0.
\end{equation}
By Duhamel's principle, the solution $v_-^\epsilon(t)$ 
satisfies 
$$v_{-}^\epsilon(t) = e^{-\epsilon t\partial_x^4}v_-(0)+ 
\int_0^t
e^{-\epsilon(t-t')\partial_x^4}\Phi(v_-^\epsilon(t'))\; dt'.$$
We have the inequality 
(by computing
$\max_{\xi \in \rr} \xi^j e^{-\epsilon t \xi^4} 
= c_j(\epsilon t)^{-j/4}$, with $c_0=1$),  
\begin{equation}
\|\partial_x^j e^{-\epsilon t\partial_x^4}f\|_{2} \leq
c_j (\epsilon t)^{-j/4}\|f\|_2 \quad j=0,1,2,3.
\end{equation}
Therefore, formally,
\begin{align*}
\|&v^{\epsilon}_{-}(t)\|_2 
\leq \|v_-(0)\|_2 
\\
&\quad+ 
\int_0^t 
\|e^{-\epsilon (t-t')\p_x^4} 
\left\{ \p_x^2 (a v_-^\epsilon)
	- \p_x (\p_x a v_-^\epsilon + 2 a \phi v_-^\epsilon) 
	+ (2 \p_x(a \phi) v_-^\epsilon + F_-)\right\} \|_2 \; dt'
\\
&\leq
\|v_-(0)\|_2 
+ 
c \int_0^t
\Big\{
\Big(\frac{1}{(\epsilon(t-t'))^{1/2}}+\frac{1}{(\epsilon(t-t'))^{1/4}} 
+ 1 \Big)\|v^{\epsilon}_{-}(t')\|_2
+ \|F_-\|_2 \Big\}\; dt' 
\\
&\leq
\|v_-(0)\|_2+c
\Big(\frac{T^{1/2}}{\epsilon^{1/2}}+\frac{T^{3/4}}{\epsilon^{1/4}} 
+ T \Big) \sup_{t \in [0,T]} \| v^{\epsilon}_-(t) \|_2
+ \int_0^T \| F_-\|_2 \; dt.
\end{align*}

A standard argument then 
shows the existence of a solution $v_{-}^\epsilon
\in C([0,T_{\epsilon}] :L^2(\mathbb R))$ to \eqref{vis}, with
$T_{\epsilon}\downarrow 0$ as $\epsilon \downarrow 0$.  Using 
the \it a priori \rm estimate \eqref{estimate}, which holds 
uniformly in $\epsilon>0$, we reapply the above 
local argument 
to extend the solution $v^{\epsilon}_{-}$ 
to the time interval $[0,T]$, with $T$ as in \eqref{TT}, for all
$\epsilon\in(0,1)$.  Letting $\epsilon\to 0$ in an appropriate
manner, we find the desired solution.

Since $v(x,t)=\varphi(x) u(x,t)$, both $u$ and $e^{\beta x}
u$ are in $C([0,T]:L^2(\mathbb R))$, with $u$ solving \eqref{2.1} in
$C([0,T]:H^{-2}(\mathbb R))$.
Also, notice that
$$
 w(x,t):=e^{\beta x}u(x,t)\in C([0,T]:L^2(\mathbb R))
$$
is a solution of the equation
$$
\aligned
\p_tw
&=i((\p_x-\beta)a(\p_x-\beta)w+W(x,t)w(x,t))
\\
&=i \p_x(a\p_x w)- 2i\beta a \p_xw +i(\beta^2 a -\beta \p_x a) w 
+ iWw,
\endaligned
$$
with
$
 w_{-}(x,0)=P_{-}(e^{\beta x}u(x,0))$, and
$w_{+}(x,T)=P_{+}(e^{\beta x}u(x,T)).
$

To prove the second part of Theorem \ref{T: T1}, 
we project the above equation onto the positive and negative 
frequencies, obtaining a coupled system for $w_\pm := P_\pm w$, 
from which we find the energy estimate
 \begin{equation}\label{a10}
\beta  \int_0^T\int a(x,t) (|D_x^{1/2}w_{+}|^2
+|D_x^{1/2}w_{-}|^2)\;dx dt\leq
 c(\|w_{-}(0)\|_2^2+\|w_{+}(T)\|_2^2).
 \end{equation}
 Therefore, from the hypothesis $a \geq \lambda > 0$, we 
see that
$w \in L^2([0,T]: H^{1/2}(\rr))$.

 We observe that formally
 $z(x,t)=D_x^{1/2}w(x,t)$ satisfies the equation
$$
 \p_t z=i\p_x(a\p_x z) -2i\beta a\p_x z +i\p_x[D_x^{1/2};a]\p_x w
 -2i\beta [D^{1/2}_x;a]\p_x w +\Gamma(z,w),
$$
 where $\Gamma(z,w)$ denotes a linear operator of 
\lq\lq order zero" in
 $(z,w)$.  Applying the projection operators, we obtain
\begin{equation}
\label{zeqn}
\begin{aligned}
 \p_t z_{\pm} &= i\p_x(a\p_x z_{\pm}) -2i\beta a\p_x z_{\pm}
+i\p_x[P_{\pm}; a]\p_x z-2i\beta[P_{\pm};a]\p_xz
 \\&\quad +P_{\pm}(i\p_x[D_x^{1/2};a]\p_x w
 -2i\beta [D^{1/2}_x;a]\p_x w +\Gamma(z,w)).
 \end{aligned}
\end{equation}
Noticing that $\partial_x = D_x^{1/2} H D_x^{1/2}$, where $H$ is the 
Hilbert transform ($\widehat{Hf}(\xi) := i \,\mbox{sgn}(\xi) \hat{f}(\xi)$), and
using Proposition 3.2, it follows that both 
 \begin{equation}\label{a11}
\begin{aligned}
 |\int (P_\pm [D_x^{1/2};a]\p_x w )\overline{z_{\pm}} \;dx|
 &=|\int ([D_x^{1/2};a] D_x^{1/2}H z) \overline{P_\pm z_{\pm}} \;dx|
\\
 &\leq c\|J^{\delta}\p_x a\|_q \,\| z\|_2\|z_\pm \|_2,
\end{aligned}
 \end{equation}
 \begin{equation}
 \label{a12}
 \begin{aligned}
 |\int (P_{\pm} (\p_x[D^{1/2}_x;a]\p_xw)) \overline{ z_{\pm}} \;dx|
&= |\int (D_x^{1/2}[D_x^{1/2};a]D_x^{1/2}Hz)
	\overline{D^{1/2}_xHz_{\pm}}\; dx|
\\
&\leq c\|J^{\delta}\p_x a\|_q \,\|D^{1/2}_x z\|_2\|D^{1/2}_x
z_{\pm}\|_2
\end{aligned}
 \end{equation}
 where we take $0< \delta < 1$ and $1 < q < \infty$ such that
$\delta>1/q$. 
Since we know 
that $\|z \|_{L^2_t L^2_x} = \|D_x^{1/2} w \|_{L^2_t L^2_x} \leq 
C_o$ ($C_o$ denoting a constant that depends on the data
$\|w_-(0)\|_2$ and $\|w_+(T)\|_2$), 
we have that $\| z(t) \|_{L^2_x} < \infty$ 
for a.e. $t$. 
Therefore, for every $\epsilon > 0$, 
we can find $t_0^\epsilon \in (0, \epsilon)$ and 
$t_1^\epsilon \in (T-\epsilon, T)$ such that 
$\| z(t_i^\epsilon) \|_{L^2_x} \leq C_o(\epsilon)$ for $i = 0, 1$. 
From the equations \eqref{zeqn},  
we obtain the following energy estimate for $z$: 
\[
\begin{aligned}
\beta \lambda \int_{t_0^\epsilon}^{t_1^\epsilon} 
\int | D_x^{1/2}z |^2 \; dx dt
&\leq
\beta \int_{t_0^\epsilon}^{t_1^\epsilon} \int a(x,t) 
(|D_x^{1/2} z_+|^2 + |D_x^{1/2} z_-|^2) \; dx dt 
\\
&\leq C_o(\epsilon) 
+ c \|J^\delta \p_x a \|_{L^\infty_t L^q_x}  
\int_{t_0^\epsilon}^{t_1^\epsilon} \| D^{1/2}_x z\|_2^2.
\end{aligned}
\]
By the hypothesis on the size of $\beta \lambda$, we can absorb the 
term on the right-hand side that arose from \eqref{a12} 
into the left-hand side.
This allows us to conclude that
$$
w\in C((0,T):H^{1/2}(\mathbb R)), \;\;\;
D_x w\in L^2(\mathbb R\times[t_0^\epsilon,t_1^\epsilon])
\;\mbox{ for every } \epsilon > 0.
$$
Reapplying this argument, 
it follows that 
$w=e^{\beta x}u\in C^{\infty}(\mathbb R\times(0,T))$.

\section{Appendix}

\begin{lemma}
\label{L: L3}

Let $T$ denote one of the following operators : $P_{+},\,P_{-}$, or $H$, 
the Hilbert transform.
Then for any $p \in (1,\infty)$ and any $l,\,m \in \mathbb  Z^+$ 
there exists $c=c(p;l;m)>0$ such that
\begin{equation}
\label{3.1}
 \| \p_x^l[T;\,a]\p_x^m f\|_p\leq c \|\p_x^{l+m} a\|_{\infty} \|f\|_p.
\end{equation}

\end{lemma}

\begin{proof}
\label{Proof of Lemma 3.1}
Without loss of generality we take $T=P_{+}$ and  observe that
$$
\p_x^l[P_{+};a] h = 
\sum_{j=0}^l\,c_{j,l}\,[P_{+};\p_x^ja]\p_x^{l-j} h,
$$
so it suffices to prove \eqref{3.1} in the case $l=0$. Also since
\begin{align*}
[P_{+}; a] 
\p_x^{m} f 
&=P_{+}(a \p_x^m f) - a P_{+}\p_x^m f 
=P_{+}(a  P_{-}\p_x^m f) + P_{+}(a \,P_{+}\p_x^m f) - a P_{+}\p_x^m f\\
&=P_{+}(a P_{-}\p_x^m f) - (I-P_{+})(a P_{+}\p_x^m f)
=P_{+}(a P_{-}\p^m_x f) - P_{-}(a P_{+}\p_x^m f),
\end{align*}
it suffices to show the inequality 
\begin{equation}
\label{3.2}
 \| P_{+}(a\,P_{-}\p_x^m f)\|_p \leq  c \|\p_x^{m} a\|_{\infty} \|f\|_p
\end{equation}
and the corresponding inequality for $\,P_{-}(a\,P_{+}\p_x^m f)\,$, 
the proof of which we omit as it is similar to the proof of 
\eqref{3.2}. 
As we commented earlier, an inequality 
related to that in \eqref{3.2} was proved in \cite{KePoVe1}.

To establish \eqref{3.2}, 
we will use the Littlewood-Paley decomposition, 
following the approach and the notation given in \cite{KePoVe1}.
First, we define functions $\eta$ and $\widetilde{\eta}$ 
centered at the frequencies $\pm 1$. 
Let $\eta \in C^\infty_0\, (\mathbb  R), \eta \geq 0, supp\; \eta
\subseteq \pm (1/2,2)$ with 
the condition $\sum\limits^\infty_{-\infty}\, \eta (2^{-k}\xi)
= 1$ for $\xi \ne 0$.
Let
$\widetilde\eta \in C^\infty_0 (\mathbb R)$, $\;\widetilde\eta \ge 0,\; supp\;
\widetilde\eta \subseteq \pm (1/8,8)$ with $\widetilde\eta (\xi) = 1$ 
for $\xi \in \pm [1/4,4]$.
Then, define 
the associated multiplication operators $Q_k$ and $\widetilde Q_k$ as follows: 
$\,(Q_k f)^\wedge (\xi) := \eta(2^{-k}\xi) \hat f (\xi)\,$ and
$(\widetilde Q_k f)^\wedge (\xi) := \widetilde\eta (2^{-k} \xi)
\hat f(\xi)\,$. 

Let $\,P_k f := \sum_{j \le k-3}\; Q_j f$;
therefore, $(P_k f) ^\wedge (\xi) = p(2^{-k}\xi) \hat f(\xi)$ 
with $p(0) = 1$ and $\,supp\; p \subseteq (-1/4,1/4)$.
Finally, define the cutoff function 
$\tilde p \in C^\infty_0 (\mathbb  R)$ with $\tilde p(\xi) = 1$ for
$\xi \in [-10,10]$ and let $(\widetilde P_k f)^\wedge(\xi) =
\tilde p(2^{-k}\xi) \hat f(\xi)$.  

Using that 
$(Q_k f)^{\wedge}$ is supported on $\pm (2^{k-1}, 2^{k+1})$
and that $(P_k f)^\wedge$ is supported on $(-2^{k-2}, 2^{k-2})$, we can
compute that 
$supp \; (Q_k f \, P_k g)^\wedge \subseteq \pm (2^{k-2}, 2^{k+2})$; 
therefore,
\begin{equation}
\label{suppid1}
Q_k f\, P_k g = \widetilde Q_k(Q_k f\, P_k g).
\end{equation}
Also, since $\widetilde P_k f = f$ if 
$supp\; \hat{f} \subset (-10 \cdot 2^k, 10 \cdot 2^k)$, 
we see that for $|j| \leq 2$,
\begin{equation}
\label{suppid2}
Q_k f\; Q_{k-j} g = \widetilde P_k (Q_k f\, Q_{k-j} g).
\end{equation}

To prove the needed estimate \eqref{3.2}, we first take the 
dyadic decomposition
of the functions on the left-hand side and split
the double sum into three parts ($l -k \leq -3$, $l -k \geq 3$, and 
$|l-k|\leq 2$):
\begin{equation*}
\begin{aligned}
\label{3.3}
&P_{+}(a\, P_{-}\p_x^m f) 
= P_{+}\Big(\sum\limits_{k,l}\, Q_k a 
P_{-}(Q_l \p_x^m  f)
\Big)
=  P_{+}\Big( \sum\limits_k\, Q_k a\, P_{-}(P_k \p_x^m f)\Big)\;+
\\
&\;\; 
P_{+}\Big(\sum\limits_k\, P_k a\, P_{-}(Q_k \p_x^m f )\Big)
+ P_{+}\Big(\sum\limits_{|j|\le 2}\, 
\sum\limits_k\, Q_k a\, P_{-}(Q_{k-j} \p_x^m f)\Big)
=: I + II + III.
\end{aligned}
\end{equation*}
 Since for all $k\in\mathbb Z$, 
$\,supp\;(P_ka \; Q_k(P_{-}\p_x^m f))^\wedge
\subset (-\infty,0)$ it follows that $II=0$.
 To estimate $I$, we use \eqref{suppid1} to write
 $$
 \aligned
 I&=\,\sum\limits_k\, P_{+}(Q_k a\, P_k (P_- \p_x^m  f))
 =\,\sum\limits_k\, \tilde Q^{+}_k(Q_k a \, P_k (P_- \p_x^m  f))\\
 &=c\sum_k\int\int e^{ix(\xi+\mu)}\,\tilde \eta^{+}(2^{-k}(\xi+\mu))\, \eta(2^{-k}\xi)\, p(2^{-k}\mu)\mu^m
 \,\hat a(\xi) \, \chi_{\rr^-}(\mu) \, \hat f(\mu)\,d\xi d\mu\\
 &=c\sum_k\int\int e^{ix(\xi+\mu)}\,m_k(\xi,\mu)\,
 \widehat {\p_x^m a}(\xi) \,(\chi_{\rr^-}(\mu)\, \hat f(\mu))d\xi d\mu,
 \endaligned
 $$
where
$m_k(\xi,\mu):=m(2^{-k}\xi,2^{-k}\mu)$, and
$m(\xi,\mu):=\tilde \eta^{+}(\xi+\mu)\,\eta(\xi)\,p(\mu)\left( \frac{\mu}{\xi}\right)^m.$

Let $q, \,h\in C^{\infty}_0(\mathbb R)$ with $q\equiv 1$ on $\,supp\;\eta$, $h\equiv 1$ on $\,supp\; p$, $\,supp\; h\subset (-1/2,1/2)$, and $\,supp \;q\subset \pm (1/4,4)$, so that
$
m(\xi,\mu)=
\tilde \eta^{+}(\xi+\mu)\,\eta(\xi)\,\mu \,p(\mu)\,\tau(\xi,\mu)$,
with
$\tau(\xi,\mu):=q(\xi) \,h(\mu)\, \mu^{m-1}/\xi^m\in C^{\infty}_0(\mathbb R^2).
$
Thus, we can write the function $\tau$ as the Fourier transform of a 
Schwartz function: 
$$
\tau(\xi,\mu)=c\int\int \,e^{i(\xi \theta +\mu \nu)}\,r(\theta,\nu)\,d\theta d\nu,\;\;\;\;\;\;\;\;\text{for some}\;\;\;\, r\in\mathbb S(\mathbb R^2).
$$
Hence,
$$
I=\int_{\nu}\int_{\theta}\,\sum_k\tilde Q_k(Q_k^{\theta}(\p_x^m a) \, P_k^{\nu}(P_-f))\;r(\theta,\nu)\,d\theta d\nu,
$$
where the symbols of $Q_k^{\theta}$ and $P_k^{\nu}$ are
$\;e^{i \theta 2^{-k}\xi}\,\eta(2^{-k}\xi)$ and 
$\,e^{i\nu 2^{-k}\mu} \,2^{-k}\mu \,p(2^{-k}\mu)$, respectively,
which belong to the class considered in \cite{KePoVe1} (page 607). 
So using Lemma A.3 in \cite{KePoVe1} 
and the Hardy-Littlewood maximal function $M$, it follows that
\begin{equation}
\label{3.4}
\begin{aligned}
\|\,\sum_k\tilde Q_k(&Q_k^{\theta}(\p_x^m a) \,P_k^{\nu}(P_-f))\|_p
\leq c
\| (\sum_k|Q_k^{\theta}(\p_x^m a) \,P_k^{\nu}(P_-f)|^2 )^{1/2}\|_p\\
& \leq c \|\sup_k|Q_k^{\theta}(\p_x^m a)|\,(\sum_k| P_k^{\nu}(P_-f)|^2)^{1/2}\|_p\\
&\leq c\|M(\p_x^m a)\|_{\infty}\;\|(\sum_k |P_k^{\nu}(P_-f)|^2)^{1/2}\|_p
\leq c\|\p_x^m a\|_{\infty}\,\|f\|_p.
\end{aligned}
\end{equation}
Finally, note that $III = 0$ if $j=-2, -1,$ or $0$.
Then, using \eqref{suppid2}, we find that
$$
III = P_{+}(\sum_{j=1}^2\,\sum_k Q_k(a)\,Q_{k-j}(P_{-}\p_x^m f))
= \sum_{j=1}^2\,\sum_k \tilde P_k^{+}(Q^*_k(\p_x^ma)\,Q^{**}_{k-j}(P_-f)),
$$
where the operators $Q_k^*$ and $Q_{k-j}^{**}$ for $j=1,2$ are given by
$$
\widehat{Q_k^*h}(\xi)
:=\frac{\eta(2^{-k}\xi)}{(2^{-k}\xi)^m}\,\hat h(\xi),
\;\;\;\;\;
\widehat{Q_{k-j}^{**}h}(\xi)
:=(2^{-k}\xi)^m\;\eta(2^{-(k-j)}\xi)\,\hat h(\xi).
$$
The symbols of these multipliers lie in the class 
considered in \cite{KePoVe1}
and $\,\tilde P_k$ is uniformly bounded in $L^p$, so an
argument similar to \eqref{3.4} provides the desired inequality.
\end{proof}

\vskip.1in

  \begin{proposition}
 \label{P: P3}
 Let  $\alpha\in[0,1),\, \beta\in(0,1)$ with $\alpha+\beta\in[0,1]$.  Then for any $p,\,q\in(1,\infty)$
 and for any $\delta>1/q$ there exists $c=c(\alpha;\beta;p;q;\delta)>0$ such that
 \begin{equation}
\label{3.10}
 \| D_x^{\alpha} [ D_x^{\beta} ; a] D_x^{1-(\alpha+\beta)} f\|_p
 \leq \| J^{\delta}\,\p_x a\|_q \|f\|_p,
 \end{equation}
where $\,J:=(1-\p_x^2)^{1/2}$.
\end{proposition}

\noindent{\it{Note.}}
The inequality \eqref{3.10} still holds with the same proof for 
$\tilde D^{s}_x = H D_x^s$ in place of $D_x^s$.
Also, 
in the case $\beta=1$, we can use 
 $
 [D_x;a]f=[H;a]\p_x f + H(\p_x a\,f)
 $
and \eqref{3.1} to obtain the inequality \eqref{3.10} 
with  $q=\infty$ and $\delta=0$.

\begin{proof}
\label{Proof of Lemma 3.2}

 We observe that
 $$
 D_x^{\alpha} [ D_x^{\beta} ; a] D_x^{1-(\alpha+\beta)} f = [D_x^{\alpha+\beta};a] D_x^{1-(\alpha+\beta)}f
 -[D_x^{\alpha};a]D_x^{1-\alpha}f.
 $$
 Therefore, it suffices to consider the case $\alpha=0$.  But the proof of this case follows by
 combining  the argument  in  Proposition A.2, Lemma A.3, and 
Theorem A.8 in the appendix of \cite{KePoVe1} with $\alpha=1$ and  
the Sobolev inequality, so it will be omitted.
\end{proof}

\bibliographystyle{amsplain}

\end{document}